# On Faster Convergence of Cyclic Block Coordinate Descent-type Methods for Strongly Convex Minimization[*]

**Xingguo Li**  lixx1661@umn.edu
*Department of Electrical and Computer Engineering*
*University of Minnesota Twin Cities*
*Minneapolis, MN 55455, USA*

**Tuo Zhao**  tourzhao@gatech.edu
*School of Industrial and Systems Engineering*
*Georgia Institute of Technology*
*Atlanta, GA 30318, USA*

**Raman Arora**  arora@cs.jhu.edu
*Department of Computer Science*
*Johns Hopkins University*
*Baltimore, MD 21218, USA*

**Han Liu**  hanliu@northwestern.edu
*Department of Electrical Engineering and Computer Science*
*Northwestern University*
*Evanston, IL 60208, USA*

**Mingyi Hong**  mhong@umn.edu
*Department of Electrical and Computer Engineering*
*University of Minnesota Twin Cities*
*Minneapolis, MN 55455, USA*

## Abstract

The cyclic block coordinate descent-type (CBCD-type) methods, which performs iterative updates for a few coordinates (a block) simultaneously throughout the procedure, have shown remarkable computational performance for solving strongly convex minimization problems. Typical applications include many popular statistical machine learning methods such as elastic-net regression, ridge penalized logistic regression, and sparse additive regression. Existing optimization literature has shown that for strongly convex minimization, the CBCD-type methods attain iteration complexity of $\mathcal{O}(p\log(1/\epsilon))$, where $\epsilon$ is a pre-specified accuracy of the objective value, and $p$ is the number of blocks. However, such iteration complexity explicitly depends on $p$, and therefore is at least $p$ times worse than the complexity $\mathcal{O}(\log(1/\epsilon))$ of gradient descent (GD) methods. To bridge this theoretical gap, we propose an improved convergence analysis for the CBCD-type methods. In particular, we first show that for a family of quadratic minimization problems, the iteration complexity $\mathcal{O}(\log^2(p) \cdot \log(1/\epsilon))$ of the CBCD-type methods matches that of the GD methods in term of dependency on $p$, up to a $\log^2 p$ factor. Thus our complexity bounds are sharper than the existing bounds by at least a factor of $p/\log^2(p)$. We also provide a lower bound to confirm that our improved complexity bounds are tight (up to a $\log^2(p)$ factor), under the assumption that the largest and smallest eigenvalues of the Hessian matrix do not scale with $p$. Finally, we generalize our analysis to other strongly convex minimization problems beyond quadratic ones.

---

[*]. Some preliminary results in this paper were presented at the 19th International Conference on Artificial Intelligence and Statistics (Li et al., 2016). This research is supported by NSF DMS1454377-CAREER; NSF IIS 1546482-BIGDATA; NIH R01MH102339; NSF IIS1408910; NSF IIS1332109; NIH R01GM083084; NSF CMMI1727757.

## 1. Introduction

We consider a class of composite convex minimization problems:

$$x^* = \underset{x \in \mathbb{R}^d}{\operatorname{argmin}} \, \mathcal{F}(x), \; \mathcal{F}(x) = \mathcal{L}(x) + \mathcal{R}(x), \tag{1}$$

where $\mathcal{L}(\cdot)$ is a twice differentiable loss function and $\mathcal{R}(\cdot)$ is a possibly nonsmooth and strongly convex penalty function. Many popular statistical machine learning problems are of the form (1), such as elastic-net regression (Zou and Hastie, 2005), ridge penalized logistic regression (Hastie et al., 2009), support vector machine (Vapnik and Vapnik, 1998) and many others (Hastie et al., 2009). For notational simplicity, we assume that there exists a partition of $d$ coordinates such that

$$x = [x_1^\top, \ldots, x_p^\top]^\top \in \mathbb{R}^d,$$

where $x_j \in \mathbb{R}^{d_j}$, $d = \sum_{j=1}^p d_j$, and $d_j \ll p$. The penalty function $\mathcal{R}(x)$ in these applications is block coordinate decomposable, i.e., $\mathcal{R}(x) = \sum_{j=1}^p \mathcal{R}_j(x_j)$. Then we can rewrite the objective in (1) as

$$\mathcal{F}(x) = \mathcal{L}(x_1, \ldots, x_p) + \sum_{j=1}^p \mathcal{R}_j(x_j).$$

Many algorithms such as gradient decent-type (GD-type) methods (Nesterov, 2004, 2007; Beck and Teboulle, 2009b,a; Becker et al., 2011), cyclic block coordinate descent-type (CBCD-type) methods (Luo and Tseng, 1992; Tseng, 1993, 2001; Friedman et al., 2007; Liu et al., 2009; Tseng and Yun, 2009; Saha and Tewari, 2013; Nutini et al., 2015; Zhao and Liu, 2015; Zhao et al., 2014b,a, 2012; Li et al., 2015b), and alternating direction method of multipliers (ADMM) (Gabay and Mercier, 1976; Boyd et al., 2011; He and Yuan, 2015; Hong and Luo, 2012; Zhao and Liu, 2012; Liu et al., 2014, 2015; Li et al., 2015a)) have been proposed to solve (1). Among these algorithms, the CBCD-type methods have been immensely successful (Friedman et al., 2007, 2010; Mazumder et al., 2011; Tibshirani et al., 2012; Razaviyayn et al., 2013; Zhao et al., 2014a). One popular instance of the CBCD-type methods is the cyclic block coordinate minimization (CBCM) method, which minimizes (1) with respect to a single block of variables while holding the rest fixed. Particularly, at the $(t+1)$-th iteration, given $x^{(t)}$, we choose to solve a collection of optimization problems: For $j = 1, \ldots, p$,

$$x_j^{(t+1)} = \underset{x_j}{\operatorname{argmin}} \, \mathcal{L}\left(x_{1:(j-1)}^{(t+1)}, x_j, x_{(j+1):p}^{(t)}\right) + \mathcal{R}_j(x_j), \tag{2}$$

where $x_{1:(j-1)}^{(t+1)}$ and $x_{(j+1):p}^{(t)}$ are defined as

$$x_{1:(j-1)}^{(t+1)} = [x_1^{(t+1)\top}, \ldots, x_{j-1}^{(t+1)\top}]^\top \text{ and } x_{(j+1):p}^{(t)} = [x_{j+1}^{(t)\top}, \ldots, x_p^{(t)\top}]^\top.$$

For some applications (e.g. elastic-net penalized linear regression), we can obtain a simple closed form solution to (2), but for many other applications (e.g. ridge-penalized logistic regression), (2) does not admit a closed form solution and requires more sophisticated optimization procedures.

A popular alternative to CBCD-type methods is to solve a quadratic approximation of (2) using the cyclic block coordinate gradient descent (CBCGD) method. For notational simplicity, we denote the partial gradient $\nabla_{x_j}\mathcal{L}(x)$ by $\nabla_j\mathcal{L}(x)$. Then the CBCGD method solves a collection of optimization problems: For $j = 1, \ldots, p$,

$$x_j^{(t+1)} = \underset{x_j}{\operatorname{argmin}}\ (x_j - x_j^{(t)})^\top \nabla_j \mathcal{L}\left(x_{1:(j-1)}^{(t+1)}, x_{j:p}^{(t)}\right) + \frac{\eta_j}{2}\|x_j - x_j^{(t)}\|^2 + \mathcal{R}_j(x_j), \qquad (3)$$

where $\eta_j > 0$ is a step-size parameter for the $j$-th block.

There have been many results on iteration complexity of block coordinate descent-type (BCD-type) methods, but most of them focus on the randomized BCD-type methods, where blocks are randomly chosen with replacement in each iteration (Shalev-Shwartz and Tewari, 2011; Richtárik and Takáč, 2012; Lu and Xiao, 2015), which demonstrate better iteration complexities than cyclic BCD-type methods in the worst case scenarios (Lee and Wright, 2016; Sun and Ye, 2016). In contrast, existing literature on cyclic BCD-type methods is rather limited. Beck and Tetruashvili (2013) focus on minimizing smooth objective functions, and has shown that given a pre-specified accuracy $\epsilon$ for the objective value, the CBCGD method attains linear iteration complexity of $\mathcal{O}(\log(1/\epsilon))$ for minimizing smooth and strongly convex problems, and sublinear iteration complexity of $\mathcal{O}(1/\epsilon)$ for smooth and nonstrongly convex problems. Hong et al. (2017); Yun (2014); Sun and Hong (2015) focus on minimizing nonsmooth composite objective functions such as (1), and has shown that the CBCM and CBCGD methods attain sublinear iteration complexity of $\mathcal{O}(1/\epsilon)$, when the objective function is nonstronlgy convex.

Here, we are interested in establishing an improved iteration complexity of the CBCM and CBCGD methods, when the nonsmooth composite objective function is strongly convex. Beck and Tetruashvili (2013) has shown that for smooth minimization, the CBCGD method attains linear iteration complexity of

$$\mathcal{O}\left(\frac{L_{\max} \cdot pL^2 \log(1/\epsilon)}{L_{\min}^2 \cdot \mu}\right), \qquad (4)$$

where $L$ is the Lipschitz constant of the gradient of the objective function, $\mu$ is the strongly convex constant of the objective function, $L_{\max} = \max_i L_i$, $L_{\min} = \min_i L_i$, and $L_i$ is the Lipschitz constant of $i$-th block of the gradient of the objective function. However, such an iteration complexity depends on $p$ (the number of blocks), and therefore is at least $p$ times worse than the complexity $\mathcal{O}\left(\mu^{-1}L\log(1/\epsilon)\right)$ of the gradient descent (GD) methods.

To bridge this theoretical gap, we propose an improved convergence analysis for the CBCD-type methods. Specifically, we show that for a family of quadratic minimization problems, the iteration complexity of the CBCD-type methods matches that of the GD methods in term of dependency on $p$ up to a $\log^2(p)$ factor. More precisely, when $\mathcal{L}(x)$ is quadratic, the iteration complexity of the CBCGD method is

$$\mathcal{O}\left(\frac{\log^2(p) L^2 \log(1/\epsilon)}{L_{\min}^\mu \cdot \mu}\right), \tag{5}$$

where $L_{\min}^\mu = \min_i\{L_i + \mu_i\}$ and $\mu_i$ is the strongly convex constant with respect to the $i$-th block of variables. Note that $L_{\min}^\mu \geq L_{\min}$. As can be seen easily, (5) is better than (4) by a factor of at least $\frac{L_{\max} \cdot p}{L_{\min} \cdot \log^2(p)}$. Note that We also provide a lower bound analysis that confirms that our improved iteration complexity is tight up to a $\log^2(p)$ factor if the largest and smallest eigenvalues of the Hessian matrix do not scale with $p$. Similar results hold for the CBCM method. We remark here that when the problem is quadratic (e.g., ridge-penalized linear regression with squared loss), (1) can be written as solving a linear system. Then the CBCM method is equivalent to the Gauss-Siedel method, which also has a linear convergence rate (Golub and Van Loan, 2012). Nevertheless, our major effort is to improve the dependence of the constant factor on the problem parameters (e.g., the block size $p$ and Lipchitz constant $L$) in the iteration complexity, which is more difficult to analyze in the blockwise minimization case for the Gauss-Siedel method[1]. In addition, the Gauss-Siedel method is not applicable beyond the quadratic case in general.

Finally, we generalize our analysis to other strongly convex minimization problems beyond quadratic minimization. Specifically, for smooth minimization, the iteration complexity of the CBCGD method is

$$\mathcal{O}\left(\frac{L_{\max}^\beta \cdot p L^\beta \log(1/\epsilon)}{L_{\min}^\beta \cdot \mu}\right), \tag{6}$$

where $L^\beta = L + \beta$, $L_{\max}^\beta = \max_i\{L_i + \beta_i\}$, $L_{\min}^\beta = \min_i\{L_i + \beta_i\}$, $\beta$ is the Lipchitz constant of gradient $\nabla \mathcal{R}(\cdot)$, and $\beta_i$ is the Lipchitz constant of $i$-th block of gradient $\nabla_i \mathcal{R}(\cdot)$ for smooth $\mathcal{R}(\cdot)$. Note that $L^\beta$, $L_{\max}^\beta$, and $L_{\min}^\beta$ are the Lipchitz constants of the gradient (and the corresponding blocks) of the objective function, which are identical to those considered in Beck and Tetruashvili (2013) when the objective is of the composite form as in (1). This indicates that (6) is better than (4) by a factor of $L^\beta / L_{\min}^\beta$, which is at least of order $\sqrt{p}$ and can be much more significant for ill-conditioned problems. Similar results hold for nonsmooth regularized minimization and their counter parts for the CBCM method; for more details refer to Table 1. It is worth mentioning that all the above results on the CBCD-type methods can be used to establish the iteration complexity for the popular permuted BCM (PBCM) and permuted BCGD (PBCGD) methods, in which the blocks are randomly sampled without replacement in each round. Improvement in terms of the dependence on constants are provided in Sun and Hong (2015) for the nonstronlgy convex problems.

The rest of the paper is organized as follows. In Section 2, we introduce notations and preliminary assumptions. Then we provide the main results of improved convergence

---

1. It requires to find an upper bound of the contraction constant $\|H_L^{-1}(H - H_L)\|$ in the Gauss-Siedel method, where $H$ is the coefficient matrix of the linear system and $H_L$ is the lower triangular matrix of $H$. Note that it requires $\|H_L^{-1}(H - H_L)\| < 1$ for the convergence of the Gauss-Siedel method.

Table 1: Compared with Beck and Tetruashvili (2013), our contributions are manyfold: (1) Developing the iteration complexity bounds of the CBCM and CBCGD methods for different specifications on $\mathcal{L}(\cdot)$ and $\mathcal{R}(\cdot)$; (2) Developing the iteration complexity bound of CBCGD for quadratic $\mathcal{L}(\cdot)$ + nonsmooth $\mathcal{R}(\cdot)$; (3) Improving the iteration complexity bound of CBCGD for smooth $\mathcal{R}(\cdot)$.

|     | Method | $\mathcal{L}(\cdot)$ | $\mathcal{R}(\cdot)$ | Our analysis | Beck & Tetruashvili (2013) |
| --- | --- | --- | --- | --- | --- |
| [a] | CBCGD | Quadratic | Smooth | $\mathcal{O}\left(\frac{\log^2(p)L^2\log(1/\epsilon)}{L_{\min}^\mu \cdot \mu}\right)$ | $\mathcal{O}\left(\frac{L_{\max}\cdot pL^2\log(1/\epsilon)}{L_{\min}^2\cdot \mu}\right)$ |
| [b] | CBCGD | Quadratic | Nonsmooth | $\mathcal{O}\left(\frac{\log^2(p)L^2\log(1/\epsilon)}{L_{\min}^\mu \cdot \mu}\right)$ | N/A |
| [c] | CBCGD | General Convex | Smooth | $\mathcal{O}\left(\frac{L_{\max}^\beta \cdot pL^\beta \log(1/\epsilon)}{L_{\min}^\beta \cdot \mu}\right)$ | $\mathcal{O}\left(\frac{L_{\max}^\beta \cdot p(L^\beta)^2 \log(1/\epsilon)}{(L_{\min}^\beta)^2 \cdot \mu}\right)$ |
| [d] | CBCGD | General Convex | Nonsmooth | $\mathcal{O}\left(\frac{L_{\max}\cdot pL\log(1/\epsilon)}{L_{\min}^\mu \cdot \mu}\right)$ | N/A |
| [e] | CBCM | Quadratic | Smooth | $\mathcal{O}\left(\frac{\log^2(p)L^2\log(1/\epsilon)}{\mu_{\min}\cdot \mu}\right)$ | N/A |
| [f] | CBCM | Quadratic | Nonsmooth | $\mathcal{O}\left(\frac{\log^2(p)L^2\log(1/\epsilon)}{\mu_{\min}\cdot \mu}\right)$ | N/A |
| [g] | CBCM | General Convex | Smooth | $\mathcal{O}\left(\frac{L_{\max}\cdot pL\log(1/\epsilon)}{\mu_{\min}\cdot \mu}\right)$ | N/A |
| [h] | CBCM | General Convex | Nonsmooth | $\mathcal{O}\left(\frac{L_{\max}\cdot pL\log(1/\epsilon)}{\mu_{\min}\cdot \mu}\right)$ | N/A |

**Remark**: See Theorem 3 for [a] and [b]; See Theorem 4 for [e] and [f]; See Theorem 7 for [c]; See Theorem 8 for [d], [g], and [h]. When $\mathcal{R}(\cdot)$ is nonsmooth, the optimization problem is actually solved by the cyclic block coordinate proximal gradient (CBCPGD) method. For notational convenience in this paper, however, we simply call it the CBCGD method.

analysis for CBCD-type approachs in Section 3. Numerical evaluations are provided in Section 4, followed by further discussions in Section 5.

## 2. Notations and Assumptions

We start with introducing notations used in this paper. Given a vector $\boldsymbol{v} = (v_1,\ldots,v_d)^\top \in \mathbb{R}^d$, we define vector norms: $\|\boldsymbol{v}\|_1 = \sum_j |v_j|$, $\|\boldsymbol{v}\|^2 = \sum_j v_j^2$, and $\|\boldsymbol{v}\|_\infty = \max_j |v_j|$. Let $\{\mathcal{A}_1,\ldots,\mathcal{A}_p\}$ be a partition of all $d$ coordinates with $|\mathcal{A}_j| = d_j$ and $\sum_{j=1}^p d_j = d$. We use $\boldsymbol{v}_j$ to denote the subvector of $\boldsymbol{v}$ with all indices in $\mathcal{A}_j$. Given a matrix $A \in \mathbb{R}^{d\times d}$, we use $\lambda_{\max}(A)$ and $\lambda_{\min}(A)$ to denote the largest and smallest eigenvalues of $A$. We denote $\|A\|$ as the spectral norm of $A$ (i.e., the largest singular value). We denote $\otimes$ and $\odot$ as the Kronecker product and Hadamard (entrywise) product for two matrices respectively.

Before we proceed with our analysis, we introduce some assumptions on $\mathcal{L}(\cdot)$ and $\mathcal{R}(\cdot)$.

**Assumption 1** *$\mathcal{L}(\cdot)$ is convex, and its gradient mapping $\nabla\mathcal{L}(\cdot)$ is Lipschitz continuous and also blockwise Lipschitz continuous, i.e., there exist positive constants $L$ and $L_j$'s such that for any $x,x' \in \mathbb{R}^d$ and $j = 1,\ldots,p$, we have*

$$\|\nabla\mathcal{L}(x') - \nabla\mathcal{L}(x)\| \leq L\|x - x'\| \text{ and } \|\nabla_j \mathcal{L}\left(x_{1:(j-1)}, x_j', x_{(j+1):p}\right) - \nabla_j\mathcal{L}(x)\| \leq L_j\|x_j - x_j'\|.$$

*We define $L_{\max} = \max_j L_j$ and $L_{\min} = \min_j L_j$.*

**Assumption 2** $\mathcal{R}(\cdot)$ *is strongly convex and also blockwise strongly convex, i.e., there exist positive constants $\mu$ and $\mu_j$'s such that for any $x, x' \in \mathbb{R}^d$ and $j = 1, \ldots, p$, we have*

$$\mathcal{R}(x) \geq \mathcal{R}(x') + (x - x')^\top \xi' + \frac{\mu}{2}\|x - x'\|^2 \quad and$$

$$\mathcal{R}_j(x_j) \geq \mathcal{R}_j(x'_j) + (x_j - x'_j)^\top \xi'_j + \frac{\mu_j}{2}\|x_j - x'_j\|^2,$$

*for all $\xi' \in \partial \mathcal{R}(x')$. We define $\mu_{\min} = \min_j \mu_j$.*

For notational simplicity, we define auxiliary variables

$$L^\mu_{\min} = \min_j L_j + \mu_j \quad \text{and} \quad y^{(t,j)} = [x^{(t)\top}_{1:(j-1)}, x^{(t-1)\top}_{j:p}]^\top, \ j = 1, \ldots, p. \tag{7}$$

We remark that $y^{(t,j)}$ serves as an intermediate variable in our analysis, which has the first $j-1$ blocks of variables updated and the remaining blocks unchanged at $t$-th iteration. Our analysis considers $L_{\min}$, $L_{\max}$, $L^\mu_{\min}$, $\mu_{\min}$, $\mu$, and $d_{\max} = \max_j d_j$ as constants, which do not scale with the block size $p$ as in existing literature (Beck and Tetruashvili, 2013).

## 3. Improved Convergence Analysis

Our analysis consists of the following three steps:

(1) Characterize the successive descent after each CBCD iteration;

(2) Characterize the gap towards the optimal objective value after each CBCD iteration;

(3) Combine (1) and (2) to establish the iteration complexity bound.

We present our analysis under different specifications on $\mathcal{L}(\cdot)$ and $\mathcal{R}(\cdot)$.

### 3.1 Quadratic Minimization

We first consider a scenario, where $\mathcal{L}(\cdot)$ is a quadratic function. Particularly, we solve

$$x^* = \operatorname*{argmin}_{x \in \mathbb{R}^d} \mathcal{L}(x) + \mathcal{R}(x) = \operatorname*{argmin}_{\substack{x_j \in \mathbb{R}^{d_j} \\ j=1,\ldots,p}} \frac{1}{2}\left\|\sum_{j=1}^p A_{*j} x_j - b\right\|^2 + \sum_{j=1}^p \mathcal{R}_j(x_j), \tag{8}$$

where $A_{*j} \in \mathbb{R}^{n \times d_j}$ for $j = 1, \ldots, p$. Typical applications of (8) in statistical machine learning include ridge regression, elastic-net penalized regression, and sparse additive regression.

We first characterize the successive descent of the CBCGD method.

**Lemma 1** *Recall that $\mathcal{F}$ is the objective defined in (1). Suppose that Assumptions 1 and 2 hold. We choose $\eta_j = L_j$ for the CBCGD method. Then for all $t \geq 1$, we have*

$$\mathcal{F}(x^{(t)}) - \mathcal{F}(x^{(t+1)}) \geq \frac{L^\mu_{\min}}{2}\|x^{(t)} - x^{(t+1)}\|^2.$$

**Proof** At $t$-th iteration, there exists a $\xi_j^{(t+1)} \in \partial \mathcal{R}_j(x_j^{(t+1)})$ satisfying the optimality condition for each sub-problem:

$$\nabla_j \mathcal{L}(y^{(t+1,j)}) + \eta_j(x_j^{(t+1)} - x_j^{(t)}) + \xi_j^{(t+1)} = 0. \tag{9}$$

Then by definition of CBCGD given in (3), we have

$$\mathcal{F}(y^{(t+1,j+1)}) \leq \mathcal{L}(y^{(t+1,j)}) + (y^{(t+1,j+1)} - y^{(t+1,j)})^\top \nabla \mathcal{L}(y^{(t+1,j)}) \\ + \frac{L_j}{2}\|y^{(t+1,j)} - y^{(t+1,j+1)}\|^2 + \mathcal{R}(y^{(t+1,j+1)}).$$

This further implies

$$\mathcal{F}(y^{(t+1,j)}) - \mathcal{F}(y^{(t+1,j+1)}) = \mathcal{L}(y^{(t+1,j)}) + \mathcal{R}(y^{(t+1,j)}) \\ \geq (y^{(t+1,j)} - y^{(t+1,j+1)})^\top \nabla \mathcal{L}(y^{(t+1,j)}) - \frac{L_j}{2}\|y^{(t+1,j)} - y^{(t+1,j+1)}\|^2 \\ + \mathcal{R}(y^{(t+1,j)}) - \mathcal{R}(y^{(t+1,j+1)}) \\ = (x_j^{(t)} - x_j^{(t+1)})^\top \nabla_j \mathcal{L}(y^{(t+1,j+1)}) - \frac{L_j}{2}\|x_j^{(t+1)} - x_j^{(t)}\|^2 + \mathcal{R}_j(x_j^{(t)}) - \mathcal{R}_j(x_j^{(t+1)}). \tag{10}$$

By Assumptions 2, we have

$$\mathcal{R}_j(x_j^{(t)}) - \mathcal{R}_j(x_j^{(t+1)}) \geq (x_j^{(t)} - x_j^{(t+1)})^\top \xi_j^{(t+1)} + \frac{\mu_j}{2}\|x_j^{(t)} - x_j^{(t+1)}\|^2. \tag{11}$$

Combining (9), (10), and (11), we have

$$\mathcal{F}(y^{(t+1,j)}) - \mathcal{F}(y^{(t+1,j+1)}) \geq \frac{L_j + \mu_j}{2}\|x_j^{(t)} - x_j^{(t+1)}\|^2. \tag{12}$$

We complete the proof via summation of (12) over $j = 1, \ldots, p$ and the definition of $L_{\min}^\mu$. ∎

Next, we characterize the gap towards the optimal objective value.

**Lemma 2** *Suppose that Assumptions 1 and 2 hold with $d \geq 2$. Then for all $t \geq 1$, we have*

$$\mathcal{F}(x^{(t+1)}) - \mathcal{F}(x^*) \leq \frac{8L^2 \log^2(3pd_{\max})}{\mu}\|x^{(t+1)} - x^{(t)}\|^2.$$

**Proof** Since $\mathcal{L}(x)$ is quadratic, its second order Taylor expansion is tight, i.e.

$$\mathcal{L}(x^*) = \mathcal{L}(x^{(t+1)}) + \langle \nabla \mathcal{L}(x^{(t+1)}), x^* - x^{(t+1)} \rangle + \frac{1}{2}\|A(x^{(t+1)} - x^*)\|^2, \tag{13}$$

where $A = [A_{*1}, \ldots, A_{*p}] \in \mathbb{R}^{n \times d}$.

Consider matrices $\widetilde{P}$ and $\widetilde{A}$, defined as

$$\widetilde{P} = \begin{bmatrix} L_1 & 0 & 0 & \ldots & 0 & 0 \\ 0 & L_2 & 0 & \ldots & 0 & 0 \\ \vdots & \vdots & \vdots & \ldots & \vdots & \vdots \\ 0 & 0 & & \ldots & 0 & L_p \end{bmatrix} \in \mathbb{R}^{p \times p} \text{ and } \widetilde{A} = \begin{bmatrix} A_{*1} & 0 & 0 & \ldots & 0 & 0 \\ 0 & A_{*2} & 0 & \ldots & 0 & 0 \\ \vdots & \vdots & \vdots & \ldots & \vdots & \vdots \\ 0 & 0 & & \ldots & 0 & A_{*p} \end{bmatrix} \in \mathbb{R}^{np \times d}.$$

For simplicity, we assume that $d_1 = ... = d_p = m = d/p$. For any $s \in \mathbb{Z}^+$, we define the lower triangular matrix $D_s \in \mathbb{R}^{s \times s}$ as

$$D_s = \begin{bmatrix} 1 & 0 & 0 & \ldots & 0 & 0 \\ 1 & 1 & 0 & \ldots & 0 & 0 \\ 1 & 1 & 1 & \ldots & 0 & 0 \\ \vdots & \vdots & \vdots & \ldots & \vdots & \vdots \\ 1 & 1 & 1 & \ldots & 1 & 1 \end{bmatrix}$$

By the definition of $L_j$, we have that for all $j = 1, \ldots, p$,

$$L_j \geq \lambda_{\max}(A_j^\top A_j).$$

which gives us the following inequality

$$\widetilde{P} \otimes I_m \succeq \widetilde{A}^\top \widetilde{A}. \tag{14}$$

To characterize the gap towards the optimal objective, we have

$$\mathcal{F}(x^{(t+1)}) - \mathcal{F}(x^*) + \frac{\mu}{2}\|x^{(t+1)} - x^*\|^2$$

$$\overset{(i)}{\leq} \langle \nabla \mathcal{L}(x^{(t+1)}), x^{(t+1)} - x^* \rangle + \mathcal{R}(x^{(t+1)}) - \mathcal{R}(x^*) + \frac{\mu}{2}\|x^{(t+1)} - x^*\|^2$$

$$\overset{(ii)}{\leq} \langle \nabla \mathcal{L}(x^{(t+1)}), x^{(t+1)} - x^* \rangle + \langle \xi^{(t+1)}, x^{(t+1)} - x^* \rangle$$

$$\overset{(iii)}{\leq} \sum_{j=1}^p \langle \nabla_j \mathcal{L}(x^{(t+1)}) - \nabla_j \mathcal{L}\left(x^{(t)}\right), x_j^{(t+1)} - x_j^* \rangle - \sum_{j=1}^p L_j \langle x_j^{(t+1)} - x_j^{(t)}, x_j^{(t+1)} - x_j^* \rangle$$

$$= \sum_{j=1}^p \left\langle \sum_{k \geq j} A_k(x_k^{(t+1)} - x_k^{(t)}), A_j(x_j^{(t+1)} - x_j^*) \right\rangle - (x^{(t+1)} - x^{(t)})^\top (\widetilde{P} \otimes I_m)(x^{(t+1)} - x^*)$$

$$\leq (x^{(t+1)} - x^{(t)})^\top \widetilde{A}^\top (D_p \otimes I_n) \widetilde{A}(x^{(t+1)} - x^*) - (x^{(t+1)} - x^{(t)})^\top (\widetilde{P} \otimes I_m)(x^{(t+1)} - x^*)$$

$$= (x^{(t+1)} - x^{(t)})^\top \left(\widetilde{A}^\top (D_p \otimes I_n) \widetilde{A} - \widetilde{P} \otimes I_m\right)(x^{(t+1)} - x^*)$$

$$\overset{(iv)}{=} (x^{(t+1)} - x^{(t)})^\top \left(\left(A^\top A - \widetilde{A}^\top \widetilde{A}\right) \odot D_d + \widetilde{A}^\top \widetilde{A} - \widetilde{P} \otimes I_m\right)(x^{(t+1)} - x^*),$$

where (i) is from (13) as $\|A(x^{(t+1)} - x^*)\|^2 \geq 0$, (ii) is from Assumption 2, (iii) is from the optimality condition to the subproblem associated with $x_j$,

$$\langle \nabla_j \mathcal{L}(y^{(t+1,j)}) + L_j(x_j^{(t+1)} - x_j^{(t)}) + \xi_j^{(t+1)}, x_j - x_j^{(t+1)} \rangle \geq 0 \text{ for any } x_j \in \mathbb{R}^m,$$

and (iv) comes from the fact that

$$\widetilde{A}^\top (D_p \otimes I_m) \widetilde{A} = \left(A^\top A - \widetilde{A}^\top \widetilde{A}\right) \odot D_d + \widetilde{A}^\top \widetilde{A}, \tag{15}$$

where $\odot$ denotes the Hadamard product. We rewrite (15) in this way because a tight upper bound of the spectral norm of R.H.S. is easier to be obtained than the spectral norm of L.H.S. after subtracting $\widetilde{P} \otimes I_m$.

For notational convenience, we define

$$B = \left(A^\top A - \widetilde{A}^\top \widetilde{A}\right) \odot D_d + \widetilde{A}^\top \widetilde{A} - \widetilde{P} \otimes I_m,$$

then we have

$$\mathcal{F}(x^{(t+1)}) - \mathcal{F}(x^*) \le (x^{(t+1)} - x^{(t)})^\top B (x^{(t+1)} - x^*) - \frac{\mu}{2}\|x^{(t+1)} - x^*\|^2. \tag{16}$$

Minimizing the R.H.S. of the above inequality over $x^*$, we obtain

$$-\mu(x^* - x^{(t+1)}) - B^\top (x^{(t+1)} - x^{(t)}) = 0.$$

which implies

$$x^* = -\frac{B^\top (x^{(t+1)} - x^{(t)})}{\mu} + x^{(t+1)}. \tag{17}$$

In addition, we have

$$\|B\| \overset{(i)}{\le} \left\| \left(A^\top A - \widetilde{A}^\top \widetilde{A}\right) \odot D_d \right\| + \|\widetilde{A}^\top \widetilde{A} - \widetilde{P} \otimes I_m\|$$
$$\overset{(ii)}{\le} \|A^\top A - \widetilde{A}^\top \widetilde{A}\| \left(1 + \frac{1}{\pi} + \frac{\log(d)}{\pi}\right) + \|\widetilde{A}^\top \widetilde{A} - \widetilde{P} \otimes I_m\|$$
$$\overset{(iii)}{\le} \left(\|A^\top A\| + \|\widetilde{A}^\top \widetilde{A}\|\right)\left(1 + \frac{1}{\pi} + \frac{\log(d)}{\pi}\right) + \|\widetilde{A}^\top \widetilde{A} - \widetilde{P} \otimes I_m\|$$
$$\overset{(iv)}{\le} 4\|A^\top A\| \left(1 + \frac{1}{\pi} + \frac{\log(d)}{\pi}\right) \overset{(v)}{\le} 4L \log(3pd_{\max}), \tag{18}$$

where (i) and (iii) are from the triangle inequality, (iv) is from (14) and the fact that $\|\widetilde{A}^\top \widetilde{A}\| \le \|A^\top A\|$, and (v) is from $\|A^\top A\| \le L$, $1 + \frac{1}{\pi} + \frac{\log(d)}{\pi} \le \log(3d)$ for all $d \ge 2$, and $d \le pd_{\max}$. Inequality (ii) follows from the result on the spectral norm of the triangular truncation operator in Angelos et al. (1992) (Theorem 1). More specifically, let us define

$$L_d = \max\left\{\frac{\|A \odot D_d\|}{\|A\|} : A \in \mathbb{R}^{d \times d}, A \ne 0\right\},$$

which is the largest ratio between the spectral norm of the triangular truncation of $A$ (the Hadamard product of $A$ and $D_d$) and the spectral norm of $A$. Then for any $d \ge 2$, we have from Angelos et al. (1992) that

$$\left|\frac{L_d}{\log d} - \frac{1}{\pi}\right| \le \frac{\left(1 + \frac{1}{\pi}\right)}{\log d}.$$

Plugging (17) into (16), we obtain

$$\mathcal{F}(x^{(t+1)}) - \mathcal{F}(x^*) \leq \frac{1}{2\mu}\|B(x^{(t+1)} - x^{(t)})\|^2 \stackrel{(i)}{\leq} \frac{\|B\|^2}{2\mu}\|x^{(t+1)} - x^{(t)}\|^2$$

$$\stackrel{(ii)}{\leq} \frac{8L^2 \log^2(3pd_{\max})}{\mu}\|x^{(t+1)} - x^{(t)}\|^2,$$

where (i) comes from the Cauchy-Schwarz inequality, (ii) is from (18). ∎

Using Lemmas 1 and 2, we establish the iteration complexity bound of the CBCGD method for minimizing (8) in the next theorem.

**Theorem 3** *Suppose that Assumptions 1 and 2 hold with $d \geq 2$. We choose $\eta_j = L_j$ for the CBCGD method. Given a pre-specified accuracy $\epsilon > 0$ of the objective value, we need at most*

$$\left\lceil \frac{\mu L^\mu_{\min} + 16L^2 \log^2(3pd_{\max})}{\mu L^\mu_{\min}} \log\left(\frac{\mathcal{F}(x^{(0)}) - \mathcal{F}(x^*)}{\epsilon}\right) \right\rceil$$

*iterations for the CBCGD method to ensure that $\mathcal{F}(x^{(t)}) - \mathcal{F}(x^*) \leq \epsilon$, where $L^\mu_{\min}$ is defined in (7).*

**Proof** Combining Lemmas 1 and 2, we obtain

$$\mathcal{F}(x^{(t)}) - \mathcal{F}(x^*) = [\mathcal{F}(x^{(t)}) - \mathcal{F}(x^{(t+1)})] + [\mathcal{F}(x^{(t+1)}) - \mathcal{F}(x^*)]$$

$$\geq \frac{L^\mu_{\min}}{2}\|x^{(t)} - x^{(t+1)}\|^2 + [\mathcal{F}(x^{(t+1)}) - \mathcal{F}(x^*)]$$

$$\geq \left(1 + \frac{L^\mu_{\min}\mu}{16L^2 \log^2(3pd_{\max})}\right)[\mathcal{F}(x^{(t+1)}) - \mathcal{F}(x^*)].$$

Recursively applying the above inequality for $t \geq 1$, we obtain

$$\frac{\mathcal{F}(x^{(t)}) - \mathcal{F}(x^*)}{\mathcal{F}(x^{(0)}) - \mathcal{F}(x^*)} \leq \left(1 - \frac{\mu L^\mu_{\min}}{\mu L^\mu_{\min} + 16L^2 \log^2(3pd_{\max})}\right)^t.$$

To ensure $\mathcal{F}(x^{(t)}) - \mathcal{F}(x^*) \leq \epsilon$, we only need a large enough $t$ to ensure that

$$\left(1 - \frac{\mu L^\mu_{\min}}{\mu L^\mu_{\min} + 16L^2 \log^2(3pd_{\max})}\right)^t [\mathcal{F}(x^{(0)}) - \mathcal{F}(x^*)] \leq \epsilon. \quad (19)$$

We complete the proof by combining (19) and the basic inequality $\kappa \geq \log^{-1}\left(\frac{\kappa}{\kappa-1}\right)$. ∎

As can be seen in Theorem 3, the iteration complexity depends on $p$ only in the order of $\log^2(p)$, which is generally mild in practice. The iteration complexity of the CBCM method can be established in a similar manner.

**Theorem 4** *Suppose that Assumptions 1 and 2 hold with $d \geq 2$. Given a pre-specified accuracy $\epsilon$, we need at most*

$$\left\lceil \frac{\mu\mu_{\min} + 64L^2 \log^2(3pd_{\max})}{\mu\mu_{\min}} \log\left(\frac{\mathcal{F}(x^{(0)}) - \mathcal{F}(x^*)}{\epsilon}\right) \right\rceil$$

*iterations for the CBCM method such that $\mathcal{F}(x^{(t)}) - \mathcal{F}(x^*) \leq \epsilon$*

**Proof** The overall proof also consists of three major steps: (i) successive descent, (ii) gap towards the optimal objective value, and (iii) iteration complexity.

**Successive Descent**: At $t$-th iteration, there exists a $\xi_j^{(t+1)} \in \partial \mathcal{R}_j(x_j^{(t+1)})$ satisfying the optimality condition:

$$\nabla_j \mathcal{L}(y^{(t+1,j+1)}) + \xi_j^{(t+1)} = 0. \tag{20}$$

Then we have

$$\mathcal{F}(y^{(t+1,j)}) - \mathcal{F}(y^{(t+1,j+1)}) \overset{(i)}{\geq} (x_j^{(t)} - x_j^{(t+1)})^\top \nabla_j \mathcal{L}(y^{(t+1,j+1)}) + \mathcal{R}_j(x_j^{(t)}) - \mathcal{R}_j(x_j^{(t+1)})$$

$$\overset{(ii)}{\geq} \left(\nabla_j \mathcal{L}(y^{(t+1,j+1)}) + \xi_j^{(t+1)}\right)^\top (x_j^{(t)} - x_j^{(t+1)}) + \frac{\mu_j}{2} \|x_j^{(t)} - x_j^{(t+1)}\|^2$$

$$\overset{(iii)}{=} \frac{\mu_j}{2} \|x_j^{(t)} - x_j^{(t+1)}\|^2, \tag{21}$$

where (i) is from the convexity of $\mathcal{L}(\cdot)$, (ii) is from Assumptions 2, and (iii) is from (20). By summation of (12) over $j = 1, \ldots, p$, we have

$$\mathcal{F}(x^{(t)}) - \mathcal{F}(x^{(t+1)}) \geq \frac{\mu_{\min}}{2} \|x^{(t)} - x^{(t+1)}\|^2.$$

**Gap towards the Optimal Objective Value**: The proof follows the same arguments with the proof of Lemma 2, with a few differences.

First, with the optimality condition to the subproblem associated with $x_j$,

$$\langle \nabla_j \mathcal{L}(x^{(t)}) + \xi_j^{(t+1)}, x_j - x_j^{(t+1)} \rangle \geq 0 \text{ for any } x_j \in \mathbb{R}^m,$$

we have

$$\mathcal{F}(x^{(t+1)}) - \mathcal{F}(x^*) \leq (x^{(t+1)} - x^{(t)})^\top B(x^{(t+1)} - x^*) - \frac{\mu}{2}\|x^{(t+1)} - x^*\|^2,$$

where $B = \left(A^\top A - \widetilde{A}^\top \widetilde{A}\right) \odot D_d + \widetilde{A}^\top \widetilde{A}$.

Then, using the same technique to bound the eigenvalues for matrices with Hadamard product, we have

$$\mathcal{F}(x^{(t+1)}) - \mathcal{F}(x^*) \leq \frac{L^2 \log^2(3d) + L_{\max}^2}{\mu} \|x^{(t+1)} - x^{(t)}\|^2 \leq \frac{2L^2 \log^2(3d)}{\mu} \|x^{(t+1)} - x^{(t)}\|^2.$$

**Iteration Complexity**: The analysis follows from the counter part of Theorem 3. ∎

Theorem 4 establishes that the iteration complexity of the CBCM method matches that of the CBCGD method. To the best of our knowledge, Theorems 3 and 4 are the sharpest iteration complexity analysis of the CBCD-type methods for minimizing (8). We further provide an example to establish the tightness of the above result in Appendix A.

## 3.2 General Smooth Minimization

We next consider general strongly convex smooth minimization, which includes Beck and Tetruashvili (2013) as a special case with $\mathcal{R}(x) = 0$. Here we require $\mathcal{R}(x)$ to be smooth and strongly convex.

**Assumption 3** $\mathcal{R}(\cdot)$ *is smooth and also blockwise smooth, i.e., there exist positive constants* $\beta$ *and* $\beta_j$*'s such that for* $x, x' \in \mathbb{R}^d$ *and* $j = 1, \ldots, p$*, we have*

$$\mathcal{R}(x) \leq \mathcal{R}(x') + (x - x')^\top \nabla \mathcal{R}(x') + \frac{\beta}{2}\|x - x'\|^2 \text{ and}$$

$$\mathcal{R}_j(x_j) \leq \mathcal{R}_j(x'_j) + (x_j - x'_j)^\top \nabla_j \mathcal{R}(x') + \frac{\beta_j}{2}\|x_j - x'_j\|^2.$$

*Moreover, we define* $\beta_{\max} = \max_j \beta_j$.

Moreover, we assume that the Hessian matrix $H$ of the objective function $\mathcal{F}$ exists, which is denoted as $H_{ij}(x) = \frac{\partial \mathcal{F}(x)}{\partial x_i \partial x_j}$.

Since the objective function is globally smooth, the CBCGD method can directly take the update form: For $j = 1, \ldots, p$,

$$x_j^{(t+1)} = x_j^{(t)} - \eta_j \left( \nabla_j \mathcal{L}(y^{(t+1,j+1)}) + \nabla \mathcal{R}_j(x_j^{(t)}) \right),$$

where $\eta_j > 0$ is a step-size parameter for the $j$-th block.

Typical applications of the general strongly convex smooth minimization in statistical machine learning includes ridge penalized logistic regression, and ridge penalized multinomial regression. It is worth mentioning that our analysis for the general case is applicable to smooth quadratic minimization, but is very different from the analysis in previous sections for quadratic minimization.

We first characterize the successive descent after each coordinate gradient descent (CGD) iteration.

**Lemma 5** *Suppose that Assumptions 1 and 3 hold. We choose* $\eta_j = L_j + \beta_j$ *for the CBCGD method. Then for all* $t \geq 1$*, there exists* $z^{(t,j)}$ *in the line segment of* $(x^{(t)}, y^{(t,j)})$ *for each* $j \in \{1, \ldots, p\}$ *such that*

$$\mathcal{F}(x^{(t)}) - \mathcal{F}(x^{(t+1)}) \geq \frac{\|\nabla \mathcal{F}(x^{(t)})\|^2}{2\left(L_{\max}^\beta + \frac{\|H\|^2}{L_{\min}^\beta}\right)},$$

*where* $H$ *is defined as*

$$H \triangleq \begin{bmatrix} 0 & 0 & 0 & \ldots & 0 & 0 \\ H_{21} & 0 & 0 & \ldots & 0 & 0 \\ H_{31} & H_{32} & 0 & \ldots & 0 & 0 \\ \vdots & \vdots & \vdots & \ldots & \ddots & \vdots \\ H_{p1} & H_{p2} & H_{p3} & \ldots & H_{p,p-1} & 0 \end{bmatrix}, \tag{22}$$

*with* $H_{ji} \triangleq H_{ji}(z^{(t,j)}) = \frac{\partial \mathcal{F}(z^{(t,j)})}{\partial z_j \partial z_i}$, *and* $L_{\min}^\beta$ *and* $L_{\max}^\beta$ *are defined as*

$$L_{\min}^\beta = \min\{L_j^\beta = L_j + \beta_j, j = 1, \ldots, p\} \quad \text{and} \quad L_{\max}^\beta = \max\{L_j^\beta = L_j + \beta_j, j = 1, \ldots, p\}.$$

**Proof** We first provide a lower bound of the successive descent using the gradient of $\mathcal{F}(\cdot)$ based on the Lipschitz continuity of $\nabla \mathcal{F}(\cdot)$. We have that $y^{(t,j)}$ and $y^{(t,k+1)}$ only differ at the $k$-th coordinate, and $\nabla_j \mathcal{F}(y^{t,k})$ has Lipschitz gradient with Lipschitz constant $F_j$, which implied

$$\mathcal{F}(y^{(t,j+1)}) \leq \mathcal{F}(y^{(t,j)}) + (y^{(t,j+1)} - y^{(t,j)})^\top \nabla_j \mathcal{F}(y^{(t,j)}) + \frac{F_j}{2}\|y^{(t,j+1)} - y^{(t,j)}\|^2$$

$$\stackrel{(i)}{=} \mathcal{F}(y^{(t,j)}) - \frac{2L_j^\beta - F_j}{2(L_j^\beta)^2} \|\nabla_j \mathcal{F}(y^{(t,j)})\|^2 \stackrel{(ii)}{\leq} \mathcal{F}(y^{(t,j)}) - \frac{1}{2L_j^\beta}\|\nabla_j \mathcal{F}(y^{(t,j)})\|^2,$$

where (i) is from that $x_j^{(t+1)} = x_j^{(t)} - \frac{\nabla_j \mathcal{F}(y^{(t,j)})}{L_j^\beta}$, and (ii) is from the fact that $L_j^\beta \geq F_j$. Then the decrease of the objective is

$$\mathcal{F}(x^{(t)}) - \mathcal{F}(x^{(t+1)}) = \sum_{k=1}^{p} \mathcal{F}(y^{(t,j)}) - \mathcal{F}(y^{(t,j+1)}) \geq \sum_{k=1}^{p} \frac{1}{2L_j^\beta} \|\nabla_j \mathcal{F}(y^{(t,j)})\|^2. \quad (23)$$

For simplicity, we assume that $d_1 = \ldots = d_p = m = d/p$. By the Mean Value Theorem, there exists $z^{(t,j)}$ such that

$$\nabla_j \mathcal{F}(x^{(t)}) = \nabla_j \mathcal{F}(x^{(t)}) - \nabla_j \mathcal{F}(y^{(t,j)}) + \nabla_j \mathcal{F}(y^{(t,j)})$$

$$\stackrel{(i)}{=} \nabla(\nabla_j \mathcal{F}(z^{(t,j)}))^\top (x^{(t)} - y^{(t,j)}) + \nabla_j \mathcal{F}(y^{(t,j)})$$

$$= \left[\frac{\partial \mathcal{F}(z^{(t,j)})}{\partial z_j \partial z_1}, \ldots, \frac{\partial \mathcal{F}(z^{(t,j)})}{\partial z_j \partial z_{j-1}}, 0, \ldots, 0\right] \left[\frac{\left(x_1^{(t)} - x_1^{(t+1)}\right)^\top}{L_1^\beta}, \ldots, \frac{\left(x_{j-1}^{(t)} - x_{j-1}^{(t+1)}\right)^\top}{L_{j-1}^\beta}, 0, \ldots, 0\right]^\top + \nabla_j \mathcal{F}(y^{(t,j)})$$

$$= \left[\frac{H_{j1}}{\sqrt{L_1^\beta}}, \ldots, \frac{H_{j,j-1}}{\sqrt{L_{j-1}^\beta}}, 0, \ldots, 0\right] \left[\frac{\left(x_1^{(t)} - x_1^{(t+1)}\right)^\top}{\sqrt{L_1^\beta}}, \ldots, \frac{\left(x_{j-1}^{(t)} - x_{j-1}^{(t+1)}\right)^\top}{\sqrt{L_{j-1}^\beta}}, 0, \ldots, 0\right]^\top + \nabla_j \mathcal{F}(y^{(t,j)})$$

$$= \left[\frac{H_{j1}}{\sqrt{L_1^\beta}}, \ldots, \frac{H_{j,j-1}}{\sqrt{L_{j-1}^\beta}}, \sqrt{L_j^\beta} \cdot I_m, 0, \ldots, 0\right] \left[\frac{\nabla_1 \mathcal{F}(y^{(t,1)})^\top}{\sqrt{L_1^\beta}}, \ldots, \frac{\nabla_p \mathcal{F}(y^{(t,p)})^\top}{\sqrt{L_p^\beta}}\right]^\top = h_j^\top f$$

where (i) is from the mean-value theorem, $h_j = \left[\frac{H_{j1}}{\sqrt{L_1^\beta}}, \ldots, \frac{H_{j,j-1}}{\sqrt{L_{j-1}^\beta}}, \sqrt{L_j^\beta} \cdot I_m, 0, \ldots, 0\right]^\top$ and $f = \left[\frac{\nabla_1 \mathcal{F}(y^{(t,1)})^\top}{\sqrt{L_1^\beta}}, \ldots, \frac{\nabla_p \mathcal{F}(y^{(t,p)})^\top}{\sqrt{L_p^\beta}}\right]^\top$. Let $\widetilde{H}$ be

$$\widetilde{H} = \begin{bmatrix} h_1^\top \\ \vdots \\ h_p^\top \end{bmatrix} = \begin{bmatrix} \sqrt{L_1^\beta} \cdot I_m & 0 & 0 & \ldots & 0 & 0 \\ \frac{H_{21}(z^{(t,2)})}{\sqrt{L_1^\beta}} & \sqrt{L_2^\beta} \cdot I_m & 0 & \ldots & 0 & 0 \\ \frac{H_{31}(z^{(t,3)})}{\sqrt{L_1^\beta}} & \frac{H_{32}(z^{(t,3)})}{\sqrt{L_2^\beta}} & \sqrt{L_3^\beta} \cdot I_m & \ldots & 0 & 0 \\ \vdots & \vdots & \vdots & \ldots & \ddots & \vdots \\ \frac{H_{p1}(z^{(t,p)})}{\sqrt{L_1^\beta}} & \frac{H_{p2}(z^{(t,p)})}{\sqrt{L_2^\beta}} & \frac{H_{p3}(z^{(t,p)})}{\sqrt{L_3^\beta}} & \ldots & \frac{H_{p,p-1}(z^{(t,p)})}{\sqrt{L_{p-1}^\beta}} & \sqrt{L_p^\beta} \cdot I_m \end{bmatrix}$$

Then we have

$$\|\nabla \mathcal{F}(x^{(t)})\|^2 = \sum_{j=1}^{p} \|\nabla_j \mathcal{F}(x^{(t)})\|^2 = \sum_{j=1}^{p} \|h_j^\top f\|^2 = \|\widetilde{H} f\|^2$$
$$\leq \|\widetilde{H}\|^2 \|f\|^2 = \|\widetilde{H}\|^2 \sum_{k=1}^{p} \frac{1}{2 L_j^\beta} \|\nabla_j \mathcal{F}(y^{(t,j)})\|^2. \quad (24)$$

Let $\widetilde{P}$ be defined as in the proof of Lemma 2. Then we have

$$\|\widetilde{H}\|^2 = \|\widetilde{P}^{1/2} + H \widetilde{P}^{-1/2}\|^2 \leq 2 \left( \|\widetilde{P}^{1/2}\|^2 + \|H \widetilde{P}^{-1/2}\|^2 \right) \leq 2 \left( L_{\max}^\beta + \frac{\|H\|^2}{L_{\min}^\beta} \right), \quad (25)$$

Combining (23), (24) and (25), we have

$$\mathcal{F}(x^{(t)}) - \mathcal{F}(x^{(t+1)}) \geq \sum_{k=1}^{p} \frac{1}{2 L_j^\beta} \|\nabla_j \mathcal{F}(y^{(t,j)})\|^2 \geq \frac{\|\nabla \mathcal{F}(x^{(t)})\|^2}{\|\widetilde{H}\|^2} \geq \frac{\|\nabla \mathcal{F}(x^{(t)})\|^2}{2 \left( L_{\max}^\beta + \frac{\|H\|^2}{L_{\min}^\beta} \right)}.$$

∎

We now characterize the gap towards the optimal objective after each CGD iteration.

**Lemma 6** *Suppose that Assumptions 1 and 2 hold. Then, for all $t \geq 1$, we have*

$$\mathcal{F}(x^{(t)}) - \mathcal{F}(x^*) \leq \frac{\|\nabla \mathcal{F}(x^{(t)})\|^2}{2\mu}.$$

**Proof** From the convexity of $\mathcal{L}(\cdot)$ and strong convexity of $\mathcal{R}(\cdot)$, we have

$$\mathcal{F}(x^{(t)}) - \mathcal{F}(x^*) \leq (x^{(t)} - x^*)^\top \nabla \mathcal{F}(x^{(t)}) - \frac{\mu}{2} \|x^{(t)} - x^*\|^2.$$

Minimizing the right hand side over $x^*$, we have $x^* = x^{(t)} - \frac{\nabla \mathcal{F}(x^{(t)})}{\mu}$ and the desired result. ∎

Combining the two lemmas above, we establish the iteration complexity bound of CGD.

**Theorem 7** *Suppose that Assumption 1, 2 and 3 hold. We choose $\eta_j = L_j + \beta_j$. Then, given a pre-specified accuracy $\epsilon$, we need at most*

$$\left\lceil \frac{L_{\max}^\beta + \frac{p L^\beta L_{\max}^\beta}{L_{\min}^\beta}}{\mu} \log \left( \frac{\mathcal{F}(x^{(0)}) - \mathcal{F}(x^*)}{\epsilon} \right) \right\rceil$$

*iterations such that $\mathcal{F}(x^{(t)}) - \mathcal{F}(x^*) \leq \epsilon$, where $L^\beta = L + \beta$.*

**Proof** We first bound $\|H\|$. Since $\mathcal{F}$ is convex, we have that $\nabla^2 \mathcal{F}$ is positive semi-definite (PSD). This implies that there exists a matrix $A$ such that $\nabla^2 \mathcal{F} = AA^\top$, where $A$ can be written as

$$A = \left[A_{1*}^\top, \ldots, A_{p*}^\top\right]^\top$$

and $A_{i*}$ is the $i$-th row submatrix of $A$ with $\|A_{i*}\| \leq \sqrt{L_i^\beta}$ and $\|A\| \leq \sqrt{L^\beta}$. Then we have

$$\|H\|^2 = \|H^\top H\| \overset{(i)}{\leq} \sum_{i=1}^p \|(H^\top H)_{ii}\| = \sum_{i=1}^p \|H_{*i}\|^2 \leq \sum_{i=1}^p \|\nabla_{*i}^2 \mathcal{F}\|^2 = \sum_{i=1}^p \|A_{i*}A\|^2$$

$$\leq \sum_{j=1}^p \|A_{*i}\|^2 \|A\|^2 \leq pL^\beta L_i^\beta \leq pL^\beta L_{\max}^\beta, \tag{26}$$

where (i) is from the norm compression inequality for block partitioned PSD matrix Horn and Johnson (2012) (Section 3.5). Thus we only need to combine Lemmas 5 and 6, and complete the proof by following similar lines to the proof of Theorem 3. ∎

As can be seen from Theorem 7, the established iteration complexity bound is sharper than that in Beck and Tetruashvili (2013) by a factor of $L^\beta/L_{\min}^\beta$, which is at least of order $\sqrt{p}$ in generic settings and can be $\gg \sqrt{p}$ for ill-conditioned problems.

### 3.3 General Nonsmooth Minimization

We provide an iteration complexity bound of the CBCM and CBCGD methods for a general $\mathcal{L}(\cdot)$ and a nonsmooth $\mathcal{R}(\cdot)$.

**Theorem 8** *Suppose that Assumptions 1 and 2 hold. We choose $\eta_j = L_j$ for the CBCGD method. Then given a pre-specified accuracy $\epsilon$ of the objective value, we need at most*

$$\left\lceil \frac{\mu L_{\min}^\mu + 4pL \cdot L_{\max}}{\mu L_{\min}^\mu} \log\left(\frac{\mathcal{F}(x^{(0)}) - \mathcal{F}(x^*)}{\epsilon}\right) \right\rceil$$

*iterations for the CBCGD method and at most*

$$\left\lceil \frac{\mu \mu_{\min} + pL \cdot L_{\max}}{\mu \mu_{\min}} \log\left(\frac{\mathcal{F}(x^{(0)}) - \mathcal{F}(x^*)}{\epsilon}\right) \right\rceil$$

*iterations for the CBCM method to guarantee $\mathcal{F}(x^{(t)}) - \mathcal{F}(x^*) \leq \epsilon$.*

**Proof** The three major steps are as follows.

**Successive Descent**: For CBCGD, using the same analysis of Lemma 1, we have that for all $t \geq 1$,

$$\mathcal{F}(x^{(t)}) - \mathcal{F}(x^{(t+1)}) \geq \frac{L_{\min}^\mu}{2} \|x^{(t)} - x^{(t+1)}\|^2.$$

For CBCM, using the same analysis of Theorem 4, we have that for all $t \geq 1$,

$$\mathcal{F}(x^{(t)}) - \mathcal{F}(x^{(t+1)}) \geq \frac{\mu_{\min}}{2} \|x^{(t)} - x^{(t+1)}\|^2.$$

**Gap towards the Optimal Objective Value**: By the strong convexity of $\mathcal{R}(\cdot)$, we have

$$\mathcal{F}(x) - \mathcal{F}(x^{(t+1)}) \geq \frac{\mu}{2} \|x - x^{(t+1)}\|^2 + (x - x^{(t+1)})^\top (\nabla \mathcal{L}(x^{(t+1)}) + \xi^{(t+1)}), \quad (27)$$

where $\xi_j^{(t+1)} \in \partial \mathcal{R}_j(x_j^{(t+1)})$. We then minimize both sides of (27) with respect to $x$ and obtain

$$\mathcal{F}(x^{(t+1)}) - \mathcal{F}(x^*) \leq \frac{\|\nabla \mathcal{L}(x^{(t+1)}) + \xi^{(t+1)}\|^2}{2\mu}, \quad (28)$$

For CBCGD, we have

$$\|\nabla \mathcal{L}(x^{(t+1)}) + \xi^{(t+1)}\|^2 \stackrel{(i)}{\leq} \sum_{j=1}^{p} \|\nabla_j \mathcal{L}(x^{(t+1)}) - \nabla_j \mathcal{L}(y^{(t+1,j+1)}) - L_j(x_j^{(t+1)} - x_j^{(t)})\|^2$$

$$\leq \sum_{j=1}^{p} 2\|\nabla_j \mathcal{L}(x^{(t+1)}) - \nabla_j \mathcal{L}(y^{(t+1,j+1)})\|^2 + 2L_j^2 \|x_j^{(t+1)} - y_j^{(t+1,j)}\|^2$$

$$\stackrel{(ii)}{\leq} \sum_{j=1}^{p} 2\|\nabla(\nabla_j \mathcal{L}(z))\| \cdot \|(x^{(t+1)} - y^{(t+1,j+1)})\|^2 + 2L_j^2 \|x_j^{(t+1)} - y_j^{(t+1,j)}\|^2$$

$$\stackrel{(iii)}{\leq} 4pL \cdot L_{\max} \|x^{(t+1)} - x^{(t)}\|^2, \quad (29)$$

where (i) comes from the optimality condition

$$\nabla_j \mathcal{L}(y^{(t+1,j+1)}) + L_j(x_j^{(t+1)} - x_j^{(t)}) + \xi_j^{(t+1)} = 0,$$

(ii) is from the mean-value theorem and Cauchy-Schwarz inequality, and (iii) is from the same argument as in the proof of Theorem 7. Combining (28) and (29), we have

$$\mathcal{F}(x^{(t+1)}) - \mathcal{F}(x^*) \leq \frac{2pL \cdot L_{\max} \|x^{(t+1)} - x^{(t)}\|^2}{\mu}.$$

For CBCM, we have

$$\|\nabla \mathcal{L}(x^{(t+1)}) + \xi^{(t+1)}\|^2 \stackrel{(i)}{\leq} \sum_{j=1}^{p} \|\nabla_j \mathcal{L}(x^{(t+1)}) - \nabla_j \mathcal{L}(y^{(t+1,j+1)})\|^2$$

$$\stackrel{(ii)}{\leq} pL \cdot L_{\max} \|x^{(t+1)} - x^{(t)}\|^2, \quad (30)$$

where (i) comes from the optimality condition

$$\nabla_j \mathcal{L}(y^{(t+1,j+1)}) + \xi_j^{(t+1)} = 0$$

and (ii) is from the same argument as (29). Combining (28) and (30), we have

$$\mathcal{F}(x^{(t+1)}) - \mathcal{F}(x^*) \leq \frac{pL \cdot L_{\max}\|x^{(t+1)} - x^{(t)}\|^2}{2\mu}.$$

**Iteration Complexity**: The analysis follows from the counter part of Theorem 3. ∎

Theorem 8 is a general result for the minimization of smooth loss function plus a non-smooth penalty function. In contrast, Beck and Tetruashvili (2013) only cover general smooth minimization. We also remark that the results of Theorem 7 and Theorem 8 are no better than their quadratic counterparts in Theorem 3 and Theorem 4 as $L \leq pL_{\max}$ in general.

### 3.4 Extensions to Nonstrongly Convex Minimization

For nonstrongly convex minimization, we only need to add a strongly convex perturbation to the objective function

$$\widehat{x} = \operatorname{argmin} \mathcal{F}(x) + \frac{\sigma}{2}\|x\|^2, \tag{31}$$

where $\sigma > 0$ is a perturbation parameter. Then, the results above can be used to analyze the CBCD-type methods for minimizing (31). Eventually, by setting $\sigma$ as a reasonable small value, we can establish $\mathcal{O}(1/\epsilon)$-type iteration complexity bounds up to a $\log(1/\epsilon)$ factor. See Shalev-Shwartz and Zhang (2014) for more details.

## 4. Numerical Results

We consider two typical statistical machine learning problems as examples to illustrate our analysis.

**(I) Elastic-net Penalized Linear Regression:** Let $A \in \mathbb{R}^{n \times d}$ be the design matrix, and $b \in \mathbb{R}^n$ be the response vector. We solve the following optimization problem

$$\min_{x \in \mathbb{R}^d} \quad \frac{1}{2n}\|b - Ax\|^2 + \lambda_1\|x\|^2 + \lambda_2\|x\|_1,$$

where $\lambda$ is the regularization parameter. We set $n = 10{,}000$ and $d = 20{,}000$. We simply treat each coordinate as a block (i.e., $d_{\max} = 1$). Each row of $A$ is independently sampled from a 20,000-dimensional Gaussian distribution with mean 0 and covariance matrix $\Sigma$. We randomly select 2,000 entries of $x$, each of which is independently sampled from a uniform distribution over support $(-2, +2)$. The response vector $b$ is generated by the linear model $b = Ax + \epsilon$, where $\epsilon$ is sampled from an $n$-variate Gaussian distribution $N(0, I_n)$. We set $\lambda_1 = \sqrt{\log 1/n}$ and $\lambda_2 = \sqrt{\log d/n} \approx 0.0315$. We normalize $A$ to have $\|A_{*j}\| = \sqrt{n}$ for $j = 1, .., d$, where $A_{*j}$ denotes the $j$-th column of $A$. For the BCGD method, we choose $\eta_j = 1$. For the gradient descent method, we either choose $\eta = \lambda_{\max}\left(\frac{1}{n}A^\top A\right)$, or adaptively select $\eta$ by backtracking line search.

**(II) Ridge Penalized Logistic Regression:** We solve the following optimization problem

$$\min_{x \in \mathbb{R}^d} \frac{1}{n} \sum_{i=1}^{n} \left[ \log(1 + \exp(x^\top A_{i*})) - b_i x^\top A_{i*} \right] + \lambda \|x\|^2.$$

We generate the design matrix $A$ and regression coefficient vector $x$ using the same scheme as sparse linear regression. Again we treat each coordinate as a block (i.e., $d_{\max} = 1$). The response $b = [b_1, ..., b_n]^\top$ is generated by the logistic model $b_i = $ Bernoulli$([1 + \exp(-x^\top A_{*i})]^{-1})$. We set $\lambda = \sqrt{1/n}$. For the BCGD method, we choose $\eta_j = \frac{1}{4}$. For gradient descent methods, we choose either the step size $\eta = \frac{1}{4}\lambda_{\max}\left(\frac{1}{n} A^\top A\right)$ or adaptively select $\eta$ by backtracking line search.

We evaluate the computational performance using the number of passes over $p$ blocks of coordinates (normalized iteration complexity). For the CBCGD method, we count one iteration as one pass (all $p$ blocks). For the randomized BCGD (RBCGD) method, we count $p$ iterations as one pass (since it only updates one block in each iteration). Besides the CBCGD and RBCGD methods, we also consider a variant of the CBCGD method named the permuted BCGD (PBCGD) method, which randomly permutes all indices for the $p$ blocks in each iteration. Since the RBCGD and PBCGD methods are inherently stochastic, we report the objective values averaged over 20 different runs. Moreover, for the RBCGD method, the block of coordinates is selected uniformly at random in each iteration. We consider four different settings for both elastic-net penalized linear regression and ridge penalized logistic regression based on different choices of the covariance matrix $\Sigma$ for generating the design matrix. We always choose $\Sigma_{jj} = 1$ for $j = 1, \ldots, d$, and for any $k \neq j$, we set (I) $\Sigma_{jk} = 0$; (II) $\Sigma_{jk} = 0.5$; (III) $\Sigma_{jk} = 0.75$; (IV) $\Sigma_{jk} = 0.5^{|j-k|}$. Note that the condition number of the Hessian matrix depends on $\Sigma$. Setting (I) and (IV) tend to yield well-conditioned Hessian matrices whereas Settings (II) and (III) tend to yield a badly conditioned Hessian matrix.

Figure 1 plots the gap between the objective value and the optimal as a function of number of passes for different methods. Our empirical findings can be summarized as follows: (1) All BCD-type methods attain better performance than the GD methods; (2) When the Hessian matrix is ill conditioned, i.e., in Setting (II) and (III), the CBCGD performs worse than the RBCGD and PBCGD methods, which suggests that there is a gap between cyclic and randomized BCGD. (3) When the Hessian matrix is well conditioned (e.g., in Settings (I) and (IV)), all three BCD-type methods attain good performance, and the CBCGD method slightly outperforms the PBCGD method; (4) The CBCGD method outperforms the RBCGD method in Setting (IV).

## 5. Discussions

Existing literature has established an iteration complexity of $\mathcal{O}(L \cdot \log(1/\epsilon)/\mu)$ for the gradient descent methods when solving strongly convex composite problems. However, our analysis shows that the CBCD-type methods only attains an iteration complexity of $\mathcal{O}(pLL_{\max} \cdot \log(1/\epsilon)/(L_{\min}\mu))$. Even though our analysis further shows that the iteration complexity of the CBCD-type methods can be further improved to $\mathcal{O}(\log^2(p)L^2 \cdot \log(1/\epsilon)/(L_{\min}\mu))$ for a quadratic $\mathcal{L}(\cdot)$, there still exists a gap of factor $L \log^2 p/L_{\min}$. As our numerical experiments show, however, the CBCD-type methods can actually attain a better computational performance than the gradient methods regardless of whether $\mathcal{L}(\cdot)$ is quadratic or not, thereby suggesting that perhaps there is still room for improvement in the iteration complexity analysis of the CBCD-type methods.

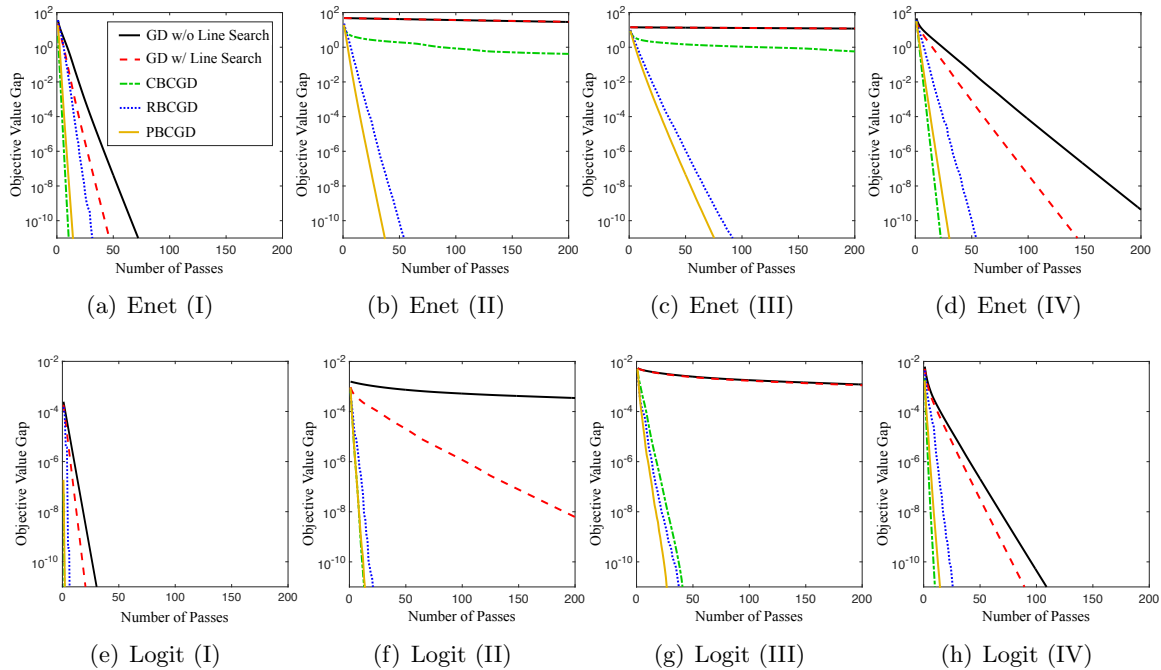

Figure 1: Comparison among different methods under different settings. "RBCGD" and "PBCGD" denote the randomized BCD-type and permuted BCD-type methods respectively. The vertical axis corresponds to the gap towards the optimal objective value, $\log[\mathcal{F}(x) - \mathcal{F}(x^*)]$; the horizontal axis corresponds to the number of passes over $p$ blocks of coordinates. Though all methods attain linear iteration complexity, their empirical behaviors are different from each others. Note that in plot ($b$) the curves for the CBCGD method and the RBCGD methods overlap.

It is also worth mentioning that though some literature claims that the CBCD-type methods works as well as the randomized BCD-type methods in practice, there do exist some counter examples, e.g. our experiment in Setting (I), where the CBCD-type methods fail significantly. This suggests that the CBCD-type methods do have some possible disadvantages in practice. To the best of our knowledge, we are not aware of any similar experimental results reported in existing literature.

Furthermore, our numerical results show that the permuted BCD-type methods, which can be viewed as a hybrid of the cyclic and the randomized BCD-type (RBCD-type) methods, has a stable performance irrespective of the problem being well conditioned or not. But to the best of our knowledge, no iteration complexity result has been established for the permuted BCD-type (PBCD-type) methods. We leave these problems for future investigation.

## Appendix A. The Tightness of the Iteration Complexity for Quadratic Problems

We next provide an example to establish the tightness of the above result. We consider the following optimization problem

$$\min_x \mathcal{H}(x) := \|Bx\|^2, \tag{32}$$

where $B \in \mathbb{R}^{p \times p}$ is a tridiagonal Toeplitz matrix defined as follows:

$$B = \begin{bmatrix} 3 & 1 & 0 & 0 & \ldots & 0 & 0 & 0 \\ 1 & 3 & 1 & 0 & \ldots & 0 & 0 & 0 \\ 0 & 1 & 3 & 1 & \ldots & 0 & 0 & 0 \\ \vdots & \vdots & \vdots & \vdots & \ddots & \vdots & \vdots & \vdots \\ 0 & 0 & 0 & 0 & \ldots & 1 & 3 & 1 \\ 0 & 0 & 0 & 0 & \ldots & 0 & 1 & 3 \end{bmatrix}.$$

Note that the minimizer to (32) is $x^* = [0, 0, \ldots, 0]^\top$, and the eigenvalues of $B$ are given by $3 + 2\cos(j\pi/(j+1))$ for $j = 1, \ldots, p$. Since the Hessian matrix of (32) is $2B^\top B$, we have

$$L = \lambda_{\max}(2B^\top B) \leq 50, \mu = \lambda_{\min}(2B^\top B) \geq 2, \mu_{\min} = 10.$$

Clearly, for this problem the largest and smallest eigenvalues of the Hessian matrix, as well as $L/\gamma$ do not scale with $p$. We consider each coordinate $x_j \in \mathbb{R}$ as a block. Then the problem can be rewritten as $\min \|\sum_{j=1}^p B_{*j} x_j\|$, where $B_{*j}$ denotes the $j$-th column of $B$. Given an initial solution $x^{(0)}$, we can show that $x^{(1)}$ is generated by

$$\begin{cases} x_1^{(1)} &= -\frac{1}{4}\left(4x_2^{(0)} + x_3^{(0)}\right), \\ x_2^{(1)} &= -\frac{1}{5}\left(4x_1^{(1)} + 4x_3^{(0)} + x_4^{(0)}\right) \\ x_3^{(1)} &= -\frac{1}{5}\left(x_1^{(1)} + 4x_2^{(1)} + x_4^{(0)} + x_5^{(0)}\right), \\ x_j^{(1)} &= -\frac{1}{5}\left(x_{j-2}^{(1)} + 4x_{j-1}^{(1)} + x_{j+1}^{(0)} + x_{j+2}^{(0)}\right), \\ x_{p-1}^{(1)} &= -\frac{1}{5}\left(x_{p-3}^{(1)} + 4x_{p-2}^{(0)} + 4x_p^{(0)}\right), \\ x_p^{(1)} &= -\frac{1}{4}\left(x_{p-2}^{(1)} + 4x_{p-1}^{(1)}\right). \end{cases} \tag{33}$$

Now we choose the initial solution

$$x^{(0)} = \left[1, \frac{9}{32}, \frac{7}{8}, 1, \ldots, 1, 1\right]^\top.$$

Then by (33), we obtain

$$x^{(1)} = \left[-\frac{1}{2}, -\frac{1}{2}, \ldots, -\frac{1}{2}, -\frac{3}{10}, -\frac{17}{40}\right]^\top,$$

which yields

$$\mathcal{H}(x^{(1)}) - \mathcal{H}(x^*) \geq \frac{25}{4}(p-3),$$

$$\|x^{(0)} - x^*\|^2 \leq p - 2 + \left(\frac{9}{32}\right)^2 + \left(\frac{7}{8}\right)^2 \leq p - 1.$$

Therefore, we have

$$\frac{\mathcal{H}(x^{(1)}) - \mathcal{H}(x^*)}{\|x^{(0)} - x^*\|^2} \geq \frac{25(p-3)}{4p} \geq \frac{22}{4}.$$

This implies that when the largest and smallest eigenvalues of the Hessian matrix do not scale with $p$ (the number of blocks), the iteration complexity is independent of $p$, and cannot be further improved. Though the independence of $p$ is only shown for the first iteration, we have similar claims in the subsequent iterations. We omit the detailed derivation due the heavy algebraic calculation.